\documentclass[12pt]{article}

\usepackage{amsmath}
\usepackage{amssymb}

\usepackage{geometry}
\usepackage{color}

\newtheorem{theorem}{Theorem}

\newtheorem{definition}{Definition}

\newtheorem{remark}{Remark}
\usepackage{mathrsfs} %

\title{{\normalsize\tt\hfill\jobname.tex}\\
On polynomial recurrence for reliability system 
with a warm reserve\\
{\small (for MPRF: in memory of A.D. Solovyev)}}
\author{A.Yu. Veretennikov\footnote{ University of Leeds, UK; National Research University Higher School of Economics, and Institute for Information Transmission Problems, Moscow, Russia, email: a.veretennikov @ leeds.ac.uk. Supported by Russian Academic Excellence Project `5-100'; funded by the RFBR grant 17-01-00633.}
}

\newif\ifabc
\abcfalse

\begin{document}

\maketitle


\begin{abstract}
Conditions for positive and  polynomial recurrence have been proposed for a class of reliability models of two elements with transitions from working state to failure and back. As a consequence, uniqueness of stationary distribution of the model is proved; the rate of convergence towards this distribution may be theoretically evaluated on the basis of the established recurrence.
\end{abstract}

\section{Introduction}
We consider a reliability model with two elements, the ``first'' and the ``second'',  each of which may be either in a working state or at a repair. 
The systems is said to be in a working state if at least one of the elements is in its working state. 
It is assumed that for each element at each state -- working or repairing -- there is an intensity of transition to another state. Independence of the elements is not assumed; instead, each intensity may depend on the states of both elements and on their elapsed times of being in the current states (working or repairing). The problem under consideration is to find conditions sufficient for the polynomial recurrence of the process, which would suffice for 
existence of a unique stationary regime and for some  bounds of the rate of  convergence to this regime. (We do not pursue the goal to establish such bounds themselves, although, in fact, recurrence to be established is rather close to it.) This problem was considered in \cite{Zv} where under certain conditions on the intensities (all of them were assumed bounded and bounded away from zero; as a matter of fact, they were also implicitly supposed to be continuous), and it was found that  the existence and uniqueness of the stationary distribution hold true, and an exponential convergence rate to it was established. This was due to an exponential moment of certain stopping times. The primary goal of this paper is to establish  {\em moment bounds} for certain stopping times under assumptions not covered by \cite{Zv}; it is known  that such bounds lead to some polynomial rates of convergence to the stationary regime for the model. 

In some earlier works for simpler models convergence was often derived via a regeneration method. As it was noted in \cite{Zv}, in this model {\em repeated} regenerations can only occur with a probability zero. However, some other technique -- called generalized regeneration -- may be used instead. Our aim is to construct Lyapunov functions which would lead to the desired a priori bounds. After this is done, we will only briefly comment about consequences for establishing convergence rates, leaving the issue to further studies. Concerning reliability theory in general, we refer to the seminal monograph \cite{GBS} and to the lecture notes \cite{Solovyev-lns} (in Russian).
The particular model consisting of two elements with constant intensities of all transitions can be found in various introductory textbooks on mathematical reliability and queueing; the case of variable  intensities belonging to an interval bounded away from zero and from infinity was treated in \cite{Zv} among other works; see the references therein.

\section{Setting}
The state space of the model is the product 
\[
{\cal S}:= \{0;1\}\times \mathbb R^+ \times\{0;1\}\times \mathbb  R^+. 
\]
The elements of ${\cal S}$ are the vectors 
$Z=(i,x;j,y)$ with $i,j = 0,1$ and $x,y \ge 0$. 
The value $i=0$ means that  the first element of the system is in the working state; the value $x$ stands for the elapsed time from the last change of the first variable $i$; the value $i=1$ signifies a failure and repairing of the first element. Similarly the values $j$ and $y$ are interpreted for the second element of the system:  $j=0$ means that  the second element of the system is in the working state; the value $y$ stands for the elapsed time from the last change of the third variable $j$; the value $j=1$ signifies a failure and repairing of the second  element.  The intensities of transitions are given by the functions $\lambda(i,x;j,y)$ for the first element and $\mu(i,x;j,y)$ for the second. It is assumed that transitions are only possible from $(i,x;j,y)$ to $(i^c,0;j,y)$, or to $(i,x;j^c,0)$, 
where by $i^c$ and $j^c$ we denote the complementary to $i$ and $j$, respectively, in the set $\{0,1\}$. 
The dependence of all intensities of the variables $i$ and $j$ are natural. Yet, their dependence of the variables $x$ and $y$ is also a frequent situation. For example, if the first element is already working for a long time ($x>\!\!>1$), then the second element may gradually (or quickly) transfer from the full or partial rest to a full readiness; then the increase of $\mu(0,x;0,y)$ in $x$ is plausible; also for $x>\!\!>1$ the values of $\mu(0,x;1,y)$ and $\mu(1,x;1,y)$ may increase in $x$ because it is likely or even surely that the second element may be required as soon as possible. Similar reasoning may be applied to the dependence of $\lambda(*,x;*,y)$ in $y$ (here $*$ signifies any value from the set $\{0,1\}$). The dependences of  $\lambda(*,x;*,y)$ in $x$ and of $\mu(*,x;*,y)$ in $y$ are even more than natural, 
for example, because of the fatigue of the elements in the working state and of the desire to finish their repairing faster if the elapsed time in the failure state becomes too long. 

By construction, $Z_t$ is a piecewise-linear Markov process continuous from the right with left limits (c\`adl\`ag) in the state space $\mathcal{S}$;  this process is also strong Markov  (see \cite{Davis}). The latter (strong Markov) property is not important in the present paper, but is rather essential in applications to the evaluation of the rate of convergence.

\subsection*{Assumptions and notations}

\noindent
{\it Suppose that there exist constants $\gamma, \Gamma>0$ such that for any $Z =(i,x;j,y)\in {\cal S}$ (state space)}
\begin{equation}\label{usl}
\begin{array}{c}
0<\frac{\gamma}{1+x} \leq \lambda(Z) \le \Gamma<\infty;
\\
\\
0<\frac{\gamma}{1+y}\leq {\mu}(Z)\le \Gamma < \infty.
\end{array}
\end{equation}
Let
\[
\Lambda(Z) := \lambda(Z) + \mu(Z).
\]
For any $Z=(i,x;j,y)$ denote (``c'' stands for a ``change'' and ``n'' for a ``no change'' for the respective variable 1 or 3):
\[
Z=((0,x),(0,y)) \, \Longrightarrow \, Z^{cn}=((1,0),(0,y)), \; Z^{nc}=((0,x),(1,0)); 
\]
\[
Z=((1,x),(0,y)) \, \Longrightarrow \, Z^{cn}=((0,0),(0,y)), \; Z^{nc}=((1,x),(1,0)); 
\]
\[
Z=((0,x),(1,y)) \, \Longrightarrow \, Z^{cn}=((1,0),(1,y)), \; Z^{nc}=((0,x),(0,0)); 
\]
\[
Z=((1,x),(1,y)) \, \Longrightarrow \, Z^{cn}=((0,0),(1,y)), \; Z^{nc}=((1,x),(0,0)).
\]
Recall that the two discrete components $i$  and $j$ in $Z=(i,x;j,y)$ cannot change both simultaneously, so that the notation like $Z^{cc}$ is not needed. 
Strictly speaking, the process $Z_t$ is not {\it regenerative} (see \cite{GZ}). More precisely, any state -- e.g., $(0,0;0,0)$ -- may be claimed the regeneration state, but the problem is that this (or any other fixed point in $\cal S$) is achievable  only with probability zero.

So, the methods of the proof of the ergodicity of the process using the theory  of regeneration processes are not directly applicable (see \cite{zv2,GZ}). Yet, extended regeneration as a base of the coupling method may be used, see \cite{Nummelin, Thorisson}. We will show how to apply the Lyapunov functions technique to this model so as to guarantee good recurrence properties eventually leading to the polynomial convergence under suitable conditions on the constants in the assumptions. Recall that the evaluation these rates is not the goal of this paper.


\subsection*{Measurable intensities and extended generator}
The standard definition and interpretation of intensities like 
\[
{\mathbb P}_{Z_t}(\mbox{exactly one jump of component $i$ on $[t,t+\Delta]$}) = \lambda(Z_t)\Delta + o(\Delta), \quad \Delta \downarrow 0, 
\]
(we stress out that the change in this event occurs just for one discrete component, not for both of them), and
\[
{\mathbb P}_{Z_t}(\mbox{exactly one jump of  component  $j$ on $[t,t+\Delta]$}) = \mu(Z_t)\Delta + o(\Delta), \quad \Delta \downarrow 0, 
\]
and also 
\[
{\mathbb P}_{Z_t}(\mbox{exactly one jump of  component $i$ or  $j$ on $[t,t+\Delta]$}) = (\lambda(Z_t)+\mu(Z_t))\Delta  + o(\Delta), 
\]
as \(\Delta \downarrow 0\)
(cf., e.g., \cite{Khintchine}) implicitly (or, in some cases explicitly) assumes that the functions $\lambda(Z)$ and $\mu(Z)$ are either constants, or, at least, continuous. However, for the discontinuous case such a definition may not be convenient if $Z_t$ happens to be the point of discontinuity of one of the functions $\lambda(\cdot)$, or  $\mu(\cdot)$.

Since we do not assume their continuity,  the definitions of intensities should be revised and reformulated more precisely. 
There are several options for that. One of them is to use a martingale approach, see, e.g., \cite[section III.5.5]{LiptserShiryaev} where it is given 
as an example 
in terms of indicators. 

\begin{definition}\label{def1}
Functions $\lambda$ and $\mu$ are called intensities of the underlying process $(Z_t,\, t\ge 0)$ iff for any smooth enough function $(h(Z), Z\in {\cal S})$ with a compact support the process
\begin{equation}\label{intdef1}
h(Z_t) - h(Z_0) - \int_0^t Lh(Z_s)\,ds 
\end{equation}
is a (possibly local) martingale, 
where $L$ is the extended generator built via the functions $\lambda$ and $\mu$ by the rule, 
\begin{align}\label{gen}
Lh(Z)
= \lambda(Z_{})(h(Z^{cn}_{})- h(Z_{})) 
+ \mu(Z)(h(Z^{nc})- h(Z))
 \nonumber \\\\ \nonumber
+ \left(\frac{\partial h}{\partial x}(Z) + \frac{\partial h}{\partial y}(Z)\right).
\end{align}
\end{definition}
Why $L$ should have this form is briefly explained below. The next slightly different but equivalent definition simultaneously highlights that the process $Z$ here is Markov, cf. \cite{Kurtz}.
\begin{definition}\label{def2}
Functions $\lambda$ and $\mu$ are called intensities of the underlying process $(Z_t,\, t\ge 0)$ iff for any smooth enough function $(h(Z), Z\in {\cal S})$ with a compact support, any positive integer $m$,  any non-random moments of time $0\le t_1 < \ldots < t_{m+1}$, and any Borel bounded functions $\phi_i$ on $\cal S$, $1\le i\le m$, 
\begin{equation}\label{intdef1}
\mathbb E\left(\left.\left(h(Z_{t_{m+1}}) - h(Z_{t_{m}}) - \int_0^t Lh(Z_s)\,ds) \,\prod_{i=1}^m \phi(Z_{t_{i}})\right) \right|{\cal F}_{t_{m}}\right) = 0,
\end{equation}
with $L$ from (\ref{gen}) and with a filtration ${\cal F}_{t}$ generated by the process $Z$ on $[0,t]$. 
\end{definition}
One more option  -- also equivalent, although, it is not our aim here to justify this equivalence  -- is to write down an explicit (although a bit cumbersome) formula for a more or less general event related to some interval of time $[0,t]$. Here we will use the convention
\[
Z+s := (i,x+s;j,y+s), \quad 
\forall \, Z = (i,x;j,y)\in {\cal S}, \; \mbox{and for any} \; s\ge 0.
\]
\begin{definition}\label{def3}
Functions $\lambda$ and $\mu$ are called intensities of the underlying process $(Z_t,\, t\ge 0)$ iff for any smooth enough function $(h(Z), Z\in {\cal S})$, any positive integer  $m$,  any $a_1, \ldots a_m$ taking valaues $1$ or $2$, any non-random moments of time $0\le s_0< s^{a_1}_1 < t^{a_1}_1 < \ldots < s^{a_m}_{m}< t^{a_m}_{m}<t$, for a generic event on $[0,t]$
\begin{align*}
A:= \{\mbox{exactly one jump of the component $i$ on each of the intervals $(s^1_*, t^1_*)$}, 
 \\
\mbox{and exactly one jump of the component $j$ on each of the intervals $(s^2_*, t^2_*)$}\}
\end{align*}
where $*$ stands for any value of the index $k=1, \ldots, m$, 
its conditional probability given ${\cal F}_{s_0}$ equals
\begin{align*}
\mathbb P(A |{\cal F}_{s_{0}}) 
= \int\limits_{s^{a_m}_m}^{t^{a_m}_m} \!\ldots\!\int\limits_{s^{a_1}_1}^{t^{a_1}_1} \prod_{k=1}^{m}\, dr^{a_k}_k
 \\\\
\times \exp\left(-\int\limits_{r_k^{a_k}}^{t_k^{a_k}} \Lambda(Z^{k}_{s^{a_k}}+\tilde r^{a_k}_k +)\,d\tilde r^{a_k}_k\right)\lambda^{a_k}_{}\left(Z^{k}_{(r^{a_k}_k)-}\right) \exp\left(-\int\limits_{s_k^{a_k}}^{r_k^{a_k}} \Lambda(Z^{k}_{s^{a_k}}+\tilde r^{a_k}_k -)\,d\tilde r^{a_k}_k\right) 
 \nonumber \\\\ \nonumber 
\times \exp\left(-\int\limits^{t}_{t_m^{a_m}} \Lambda(Z^{m}_{t^{a_m}}+ r^{}+)\,d r^{}\right)
\times \exp\left(-\int\limits_0^{s_1^{a_1}} \Lambda(Z_{0}+r_0-)\,dr_0\right), 
\end{align*}
where the filtration ${\cal F}_{t}$ is  generated by the process $Z$ on $[0,t]$, 
$\lambda^1(Z) = \lambda(Z)$, $\lambda^2(Z) = \mu(Z)$, (recall) $\Lambda(Z)=\lambda(Z)+\mu(Z)$, and the vector $Z^{k}_{r^{a_k}_k}$ for $k=1, \ldots$ is defined by induction by the rule 
\begin{align*}
Z^{k}_{(r^{a_k}_k)-} = (Z^{}_{t^{a_{k-1}}_{k-1}}+(r^{a_k}_k-t^{a_{k-1}}_{k-1})),
 \\\\
Z^{k}_{r^{a_k}_k} = (Z^{}_{t^{a_{k-1}}_{k-1}}+(r^{a_k}_k-t^{a_{k-1}}_{k-1}))^{cn}1(a_{k}=1) + (Z^{}_{t^{a_{k-1}}_{k-1}}+(r^{a_k}_k-t^{a_{k-1}}_{k-1}))^{nc}1(a_{k}=2).
\end{align*}
\end{definition}
The integration over $dr^{a_k}_k$ is performed here on the interval $(s^{a_k}_k, t^{a_k}_k)$. 
This corresponds, in particular, to the approach in  \cite{Davis}. We stress out that all the definitions lead to Dynkin's formulae below. Also, note that the following usual formulae (\ref{intu1} -- \ref{intu3b}) which are known to be valid under the assumptions of continuity of the intensities remain true without the requirement of this continuity. Here we use the convention
\[
Z+s := (i,x+s;j,y+s), \quad 
\forall \, Z = (i,x;j,y)\in {\cal S}, \; \mbox{and for any} \; s\ge 0.
\]
For any non-random values $t\ge 0, \,\Delta>0$ the following exact (not asymptotic for small $\Delta$, i.e., without $o(\Delta)$ except for (\ref{intu1a}), (\ref{intu3}) and (\ref{intu3a})) identities holds true.
\begin{equation}\label{intu1}
{\mathbb P}_{(i_t,x_t; j_t, y_t)}(\mbox{no jumps on $[t,t+\Delta]$}) = \exp\left(-\int_0^{\Delta} (\lambda+\mu)
(i_{t},x_{t}+s; j_t, y_t+s)\,ds\right).
\end{equation}
Further, 
\begin{equation}\label{intu1a}
{\mathbb P}_{(i_t,x_t; j_t, y_t)}(\mbox{more than two jumps on $[t,t+\Delta]$}) = o(\Delta).
\end{equation} 
A  complementary probability to (\ref{intu1}) is written as 
\begin{align}\label{intu2}
{\mathbb P}_{(i_t,x_t; j_t, y_t)}(\mbox{at least one jump on $[t,t+\Delta]$}) 
 \nonumber\\\\\nonumber
= 1- \exp\left(-\int_0^{\Delta} (\lambda+\mu)
(i_{t},x_{t}+s; j_t, y_t+s)\,ds\right).
\end{align}
Emphasize that both (\ref{intu1}) and (\ref{intu2}) are rigorous equalities. Respectively, 
\begin{align}\label{intu3}
{\mathbb P}_{(i_t,x_t; j_t, y_t)}(\mbox{exactly one jump on $[t,t+\Delta]$}) 
 \nonumber\\\\\nonumber
= 1- \exp\left(-\int_0^{\Delta} (\lambda+\mu)
(i_{t},x_{t}+s; j_t, y_t+s)\,ds\right) + o(\Delta);
\end{align}
and more precisely, 
\begin{align}\label{intu3a}
{\mathbb P}_{(i_t,x_t; j_t, y_t)}(\mbox{exactly one jump of component $i$ on $[t,t+\Delta]$}) 
 \nonumber\\\\\nonumber
= 1- \exp\left(-\int_0^{\Delta} \lambda
(i_{t},x_{t}+s; j_{t+s}, y_{t+s})\,ds\right) + o(\Delta),
\end{align}
and 
\begin{align}\label{intu3b}
{\mathbb P}_{(i_t,x_t; j_t, y_t)}(\mbox{at least one jump of component $i$ on $[t,t+\Delta]$}) 
 \nonumber\\\\\nonumber
= 1- \exp\left(-\int_0^{\Delta} \lambda
(i_{t},x_{t}+s; j_{t+s}, y_{t+s})\,ds\right);
\end{align}
similarly for the other component $j$, 
\begin{align}\label{intu3a}
{\mathbb P}_{(i_t,x_t; j_t, y_t)}(\mbox{exactly one jump of component $j$ on $[t,t+\Delta]$}) 
 \nonumber\\\\\nonumber
= 1- \exp\left(-\int_0^{\Delta} \mu
(i_{s},x_{s}; j_t, y_t+s)\,ds\right) + o(\Delta),
\end{align}
and 
\begin{align}\label{intu3b}
{\mathbb P}_{(i_t,x_t; j_t, y_t)}(\mbox{at least one jump of component $i$ on $[t,t+\Delta]$}) 
 \nonumber\\\\\nonumber
= 1- \exp\left(-\int_0^{\Delta} \lambda
(i_{t},x_{t}+s; j_{t+s}, y_{t+s})\,ds\right).
\end{align}

Here is a brief explanation of the form of the (extended) generator $L$ given in (\ref{gen}). In this explanation we do assume all intensities continuous. Without this assumption the formulae still can be justified, for example, by using the approach from \cite{VZ-mprf}. Let $h(Z), \, Z\in {\cal S}$ be a Borel bounded smooth enough function. For small $t>0$ and with $Z_0=Z = (i,x;j,y)$ we have by the complete expectation formula (by analogy with complete probability), 
$$
E_Zh(Z_t) = t\lambda(Z)  h(Z^{cn})
+ t\mu(Z)  h(Z^{nc}) + (1-(\lambda(Z)+ \mu(Z))t) h(Z+t) + o(t).
$$
Subtracting $h(Z)=(t\lambda(Z)+t\mu(Z))h(Z) + (1-(\lambda(Z)+ \mu(Z))t) h(Z)$ and dividing by $t$, we obtain, 
\begin{align*}
\frac{E_Zh(Z_t) - h(Z)}{t} 
= \lambda(Z)  (h(Z^{cn})-h(Z)) 
+ \mu(Z)  (h(Z^{nc}) - h(Z)) 
 \\\\
+t^{-1}(1- (\lambda(Z)+ \mu(Z))t)(h(Z+t) - h(Z))+ o(1), \quad t\downarrow 0.
\end{align*}
Since 
\begin{align*}
t^{-1}(h(Z+t) - h(Z)) \to \left(\frac{\partial h}{\partial x}(Z) + \frac{\partial h}{\partial y}(Z)\right), \quad t\downarrow 0,
\end{align*}
and 
\begin{align*}
t^{-1}(\lambda(Z)+ \mu(Z))t)(h(Z+t) - h(Z)) = o(1), \quad t\downarrow 0,
\end{align*}
we get, 
\begin{align*}
\lim_{t\downarrow 0}\frac{E_Zh(Z_t) - h(Z)}{t} 
= \lambda(Z)  (h(Z^{cn})-h(Z)) 
+ \mu(Z)  (h(Z^{nc}) - h(Z)) 
 \\\\
+ \left(\frac{\partial h}{\partial x}(Z) + \frac{\partial h}{\partial y}(Z)\right), 
\end{align*}
and this limit is uniform, as required. For the extended generator there will be no uniformity, however, we still have Dynkin's formula,
\begin{equation}\label{dyn1}
{\mathbb E}_{Z}h(Z_t) - h(Z) = 
{\mathbb E}_{Z} \int_0^t Lh(Z_{s-})\, ds,
\end{equation}
with $L$ given by (\ref{gen}).
The equation (\ref{dyn1}), ``as usual'', can be justified via the complete expectation formula; the latter with continuous intensities is a simple corollary of the convergence of Riemann's integral sums to their limit. The reader is likely to be used to the ``complete probability'' formula where the probability space is split into a no more than a countable number of events, say, $\Omega = \sum_k \Omega_k$, and then the probability of a new event $A$ equals $P(A) = \sum_k P(A\bigcap\Omega_k)$. It seems reasonable to call a similar formula for expectations by complete expectations one. Why we insist that yet for integrals some care should be taken and even that one may wish to justify such a formula accurately is that in this case $\Omega$ is split into uncountably many events, ``especially'' if the integral is Lebesgue's one. A version of such a justification of a complete expectation formula for possibly discontinuous intensities (where Lebesque's integral must be used) can be found, for example, in \cite{VZ-mprf}. 

Note that since jumps occur at each $t$ with a probability zero, the formula (\ref{dyn1}) can be rewritten in the form
\begin{equation}\label{dyn11}
{\mathbb E}_{Z}h(Z_t) - h(Z) = 
{\mathbb E}_{Z} \int_0^t Lh(Z_{s})\, ds,
\end{equation}
In turn, in terms of martingales the formula (\ref{dyn11}) (or its conditional expectation version) can be rewritten as 
\begin{equation}\label{dyn2}
h(Z_t) - h(X) - 
\int_0^t Lh(Z_s)\, ds = M_t, 
\end{equation}
with some local martingale $M_t$; if $h$ and $Lh$ are bounded, then $M_t$ in (\ref{dyn2}) is a martingale (and, in fact, this is true for a much larger class of functions $h$). 

In fact, we shall see shortly that for our purposes it is not important whether or not the martingale is local: this is because what we want to derive from it is some {\em in}equality rather than an equality. We will apply Dynkin's formula in the next sections for a justification that some function can serve as a {\em Lyapunov function}. The latter will be understood as a decrease ``on average'' along the trajectory of the process while the value of the process is not too close to the class of states $(*,0;*;0)$.

For a suitable function $h(t,Z)$ depending on $t$ and $Z$,   Dynkin's formula becomes 
\begin{equation*}\label{dyn3}
E_Zh(t,Z_t) = h(0,Z) + \int_0^t \left(Lh(s,Z_{s-})+\frac{\partial h}{\partial s}(s,Z_s)\right)\, ds,
\end{equation*}
or, equivalently (by the same reason as (\ref{dyn11})), 
\begin{equation}\label{dyn3}
E_Zh(t,Z_t) = h(0,Z) + \int_0^t \left(Lh(s,Z_s)+\frac{\partial h}{\partial s}(s,Z_s)\right)\, ds,
\end{equation}
One more equivalent version for (\ref{dyn3})  is 
\begin{equation}\label{dyn22}
h(t,Z_t) - h(0,Z_0) - 
\int_0^t \left(Lh(s,Z_s) + 
\frac{\partial h}{\partial s}(s,Z_s)
\right)\, ds = M_t, 
\end{equation}
with a (local) martingale $M_t$. If in doubt whether or not the martingale is not local, we will use some localizing sequence of stopping times in the calculus. The equations (\ref{dyn2}) and  (\ref{dyn22}) (or, more formally, their versions with $Z_{s-}$) are often called Ito's formulae.

\section{Recurrence of the process}
Let us consider the following {\em Lyapunov functions} (i.e., the functions which will be shown to possess a Lyapunov property to decrease on average outside $\mathbb K$ on the trajectory of the process), 
\[
V_{m}(Z):= (1+x +y)^m, \quad 
V_{k,m}(t,Z):= (1+t)^k (1+x +y)^m
\]
with $1\le m\le m_0$, for $Z=(i,x;j,y)$.
Further, let $K>0$, and \({\mathbb K} = {\mathbb K}(m) = {\mathbb K}(K,m) := (Z = (i,x;j,y)\in {\cal S}: \; V_{m}(Z) \le K)\), and 
\[
\tau = \tau(m) = \tau(K,m):= \inf(t\ge 0: \, Z_t \in \mathbb K).
\]

\begin{theorem}\label{Thm1}
1. (Case $k=0$,  $K>\!\!>1$) If $\gamma > 2m_0 \ge 2$, then 
\begin{equation}\label{eqtau0}
\mathbb E_{Z_0}\tau^{}(K,m_0) \le V_{m_0}(Z_0),
\end{equation}
if $K$ is large enough.
In particular, $\gamma > 2$ in (\ref{usl}) implies existence of the stationary distribution which is necessarily unique.  

\medskip

2. (Case $k > 0$,   $K>\!\!>1$) For any $k > 0$, if $\gamma>2m_0 > 2(1+2k)$, there exists a constant $C(k,K)$ such that for each $Z_0$
\begin{equation}\label{eqtauk}
\mathbb E_{Z_0}\tau^{k+1}(K,m_0) \le C(k,K)  V_{m_0}(Z_0),
\end{equation}
if $K>0$ is large enough.

\medskip

3. (Case $k > 0$,  any $K_1$) 
Under $k > 0$,  $\gamma>2m_0 > 2(1+k)$, for any $K_1$ there exists a constant $C(k,K_1)$ such that
\begin{equation}\label{eqtauk1}
\mathbb E_{Z_0}\tau^{k+1}(K_1,m_0) \le \tilde C(k_1,K) (V_{m_0}(Z)\vee (K+1)).
\end{equation}

\end{theorem}
Constants $C(k,K)$ in (\ref{eqtauk}) and $\tilde C(k,K_1)$ in (\ref{eqtauk1}) are, of course, not unique. Version of such  constant can be found below, respectively, in (\ref{ckk}) and in (\ref{est1k}). 

\begin{remark}
All values $k,m,m_0$  are  not necessarily integers. 
\end{remark} 

~

\noindent
{\bf Proof}. 
\noindent
{\bf 1.} 
Assume $m=m_0\in [1, \gamma/2)$. 
Let $N>K$, and let $T_N:= \inf(t\ge 0: \, V_{m_0}(Z_t)\ge N)$, and $\tau = \tau_K:= \inf(t\ge 0: \, Z_t\in {\mathbb K}(K))$.
Firstly let us apply Ito's formula  (\ref{dyn2}) to $V_m(Z_t)$ for $t<\tau_K \wedge T_N$. We have, 
\begin{align*}
V_m(Z_{t\wedge \tau\wedge T_N}) = V_m(Z_{0}) + \int\limits_0^{t\wedge \tau\wedge T_N}    (\lambda(Z_s)V_m(Z_s^{cn})
+\mu(Z_t)V_m(Z_s^{nc}))ds 
 \\\\
-\int\limits_0^{t\wedge \tau\wedge T_N}   ((\lambda(Z_s) +\mu(Z_s))
V_m(Z_s) 
+ 2m((1+x_s +y_s)^{m-1}))\,ds + M_{t\wedge \tau\wedge T_N}
 \\\\
= - \!\!\int\limits_0^{t\wedge \tau\wedge T_N}  \!\! (\lambda(Z_s)(V_m(Z_s) - V_m(Z_s^{cn})) -   \mu(Z_s)(V_m(Z_s) - V_m(Z_s^{nc})) 
+ 2mV_{m-1}(Z_s))\,ds 
 \\\\
+ M_{t\wedge \tau\wedge T_N},
\end{align*}
for some (possibly local) martingale $M_t$; however, at $t\wedge \tau\wedge T_N$ it is bounded, hence, with a zero expected value. We will shortly show that the term under the integral is ``strictly negative'' if the semi-norm $V_{m}(Z_t)>K$ (is large enough): this would have been clear without the positive term $\displaystyle (\partial V_m/\partial x + \partial V_m/\partial y)(Z_t)= mV_{m-1}(Z_t)$, but even when this term is present, the negative terms dominate, because in all situations the semi-norm $V_m(Z_t)$ after any change (to $Z^{cn}_t$ or to $Z_t^{nc}$) becomes less than $V_m(Z_{t-})$. As a result, we can claim that 
\begin{equation}\label{eq15}
V_m(Z_{t\wedge \tau\wedge T_N}) - V_m(Z_0)\le - C \int_0^{t\wedge \tau\wedge T_N}V_{m-1}(Z_{s})ds 
+ M_{t\wedge \tau\wedge T_N}
\end{equation}
with some $C>0$. 
Indeed, we estimate for $Z=(i,x;j,y)$, 
\begin{align*}
(V_m(Z) - V_m(Z^{cn})) 
= (1+x +y)^m  - (1+y)^m  
 \\\\
= ((1+x +y)  - (1+y))((1+x +y)^{m-1} + \ldots + (1+y)^{m-1})
 \\\\
\ge  x(1+x+y)^{m - 1} 
= xV_{m - 1}(Z), 
\end{align*}
and 
\begin{align*}
(V_m(Z) - V_m(Z^{nc})) 
= (1+x +y)^m  - (1+x)^m  
 \\\\
\ge  y (1+x+y)^{m - 1} 
= y V_{m - 1}(Z).
\end{align*}
These inequalities are trivial for any integer natural $m$ following from the ``simplified multiplication formulae''. For any $m\ge 1$ not necessarily integer they also easily follow from the equality, 
$$
(1+x +y)^m  - (1+x)^m - y(1+x+y)^{m-1} 
= (1+x+y)^{m-1} (1+x) - (1+x)^m \ge 0,
$$
as required (recall that $y\ge 0$). 
So, on $t<\tau_K$, 
\begin{align*}
-  \lambda(Z_t)(V_m(Z_t) - V_m(Z_t^{cn}))dt -   \mu(Z_t)(V_m(Z_t) - V_m(Z_t^{nc})) 
+ 2mV_{m-1}(Z_t) 
 \\\\
\le -\frac{\gamma}{1+x_t}x_t V_{m-1}(Z_t)
-\frac{\gamma}{1+y_t}y_tV_{m-1}(Z_t) + 2mV_{m-1}(Z_t).
\end{align*}
By our assumptions, 
on $t<\tau_K$ we have either $x_t\ge K/2$, or $y_t\ge K/2$ (or both). 
So, for any $\delta>0$ there exists $K(\delta)$ large enough such that  for any 
$K\ge K(\delta)$, 
\begin{equation}\label{Kdelta}
-\frac{\gamma}{1+x_t}x_t V_{m-1}(Z_t)
-\frac{\gamma}{1+y_t}y_tV_{m-1}(Z_t) \le - (1-\delta)\gamma V_{m-1}(Z_t)
\end{equation}
on $t<\tau_K$. We have, 
$(\gamma(1-\delta)-2m)V_{m-1}(Z_t) 
>0$ on this set. 
So, this leads  to the inequality 
\begin{equation}\label{e3}
\mathbb E_{Z_0}V_m(Z_{t\wedge \tau\wedge T_N}) + \mathbb  E_{Z_0}(t\wedge \tau\wedge T_N)\le V_m(Z_0),
\end{equation}
which is a  weakened version of (\ref{eq15}). 
Letting here $N\to\infty$ and $t\to\infty$, by virtue of the Fatou lemma we obtain 
\begin{equation}\label{e4}
\mathbb E_{Z_0}V_m(Z_{\tau_K}) 
+ \mathbb  E_{Z_0}\tau_K
\le V_m(Z_0).  
\end{equation}
Note that the inequalities (\ref{e3}) and (\ref{e4}) are valid for all $1\le m\le m_0$. 
In particular,  
\[
\mathbb E_{Z_0}\tau_K  \le V_m(Z_0), 
\quad \forall \, m\in [1,m_0],  
\]
which proves (\ref{eqtauk}). 

Further, because of the local mixing within any $\mathbb K(K)$ (since on $\mathbb K(K)$ all the intensities are bounded and bounded away from zero according to (\ref{usl})) and due to the Harris--Khasminskii principle, this implies existence of the stationary distribution of the Markov process $Z_t$ (\cite{HH,HH2}). Its uniqueness follows from the assumption (\ref{usl}) which implies the possibility of gluing two processes with possibly different stationary distributions. 

~

\noindent
{\bf 2.} 
Now assume $k >0$ and $m_0\ge m\ge 1$. Let us apply Ito's formula to $V_{k,m}(t,Z_t)$ for $t<\tau\wedge T_N$ and show that  under our assumptions 
\[
\mathbb E_{Z}\tau_K^{k+1} \le C(k,K) V_{m_0}(Z))
\]
with some constant $C(k,K) >0$ to be specified.
We have,  
\begin{align*}
V_{k,m}(t\wedge \tau\wedge T_N;Z_{t\wedge \tau\wedge T_N})  - V_{k,m}(0,Z_0)
 \\\\
= \int\limits_0^{t\wedge \tau\wedge T_N}(kV_{k-1,m}(s,Z_{s}) +
  \lambda(Z_t)V_{k,m}(t;Z_t^{cn}) 
+\mu(Z_s)V_{k,m}(s;Z_s^{nc})) \,ds
 \\\\
-\int\limits_0^{t\wedge \tau\wedge T_N}(\lambda(Z_t) +\mu(Z_t))
(V_{k,m}(s;Z_s) 
+ mV_{k,m-1}(s;Z_s))ds 
+ M_{t\wedge \tau\wedge T_N}
 \\\\
=  \int\limits_0^{t\wedge \tau\wedge T_N}\left(\frac{k}{1+s}V_{k,m}(s,Z_{s}) +
  \lambda(Z_t)V_{k,m}(t;Z_t^{cn}) 
+\mu(Z_s)V_{k,m}(s;Z_s^{nc})\right) \,ds 
 \\\\
-\int\limits_0^{t\wedge \tau\wedge T_N}(\lambda(Z_t) +\mu(Z_t))
(V_{k,m}(s;Z_s) 
+ mV_{k,m-1}(s;Z_s))ds 
+ M_{t\wedge \tau\wedge T_N}
 \\\\
= - \!\!\!\int\limits_0^{t\wedge \tau\wedge T_N} \!\!\!(\lambda(Z_s)(V_{k,m}(s;Z_s) - V_{k,m}(s;Z_s^{cn}))ds +   \mu(Z_s)(V_{k,m}(s;Z_s) - V_{k,m}(s;Z_s^{nc}))ds 
 \\\\
+ \int\limits_0^{t\wedge \tau\wedge T_N}\left(\frac{k}{1+s}V_{k,m}(s,Z_{s}) + mV_{k,m-1}(s;Z_s)\right) ds + M_{t\wedge \tau\wedge T_N}.
\end{align*}
Similarly to the calculus in the first step of the proof, this implies the following  inequality with $c=2\gamma-m$, 
\begin{align*}
\mathbb EV_{k,m}({t\wedge \tau\wedge T_N};Z_{t\wedge \tau\wedge T_N}) \le V_{k,m}(0,Z_0) 
 \\\\
-  \mathbb E\int_0^{t\wedge \tau\wedge T_N} \left(cV_{k,m-1}(s;Z_{s}) - \frac{k}{1+s}V_{k,m}(s,Z_{s})\right)\,ds, 
\end{align*}
and by the Fatou Lemma also 
\begin{equation}\label{e33}
\mathbb EV_{k,m}(\tau_K;Z_{\tau_K}) \le V_{k,m}(0,Z_0) - \mathbb  E\int\limits_0^{\tau_K} \left(cV_{k,m-1}(s;Z_{s}) - \frac{k}{1+s}V_{k,m}(s,Z_{s})\right)ds.
\end{equation}
Now consider the identity
\[
1 = 1(1+x_s+y_s>\epsilon (1+s)/k) + 1(1+x_s+y_s\le \epsilon (1+s)/k), 
\]
and insert this split of unity under the integral. Then for any $0< \epsilon < c = (1-\delta)\gamma - 2m$ the expression with the indicator $1(1+x_s+y_s\le \epsilon (1+s)/k)$ will be dominated by the term $cV_{k,m-1}$ and we will get with $c'=c -\epsilon$, 
\begin{align}
\mathbb EV_{k,m}(\tau_K;Z_{t\wedge \tau_K}) \le V_{k,m}(0,Z_0) 
 \nonumber \\\nonumber \\
- \mathbb  E\int\limits_0^{\tau_K} \left(c'V_{k,m-1}(s;Z_{s}) - \frac{k}{1+s}V_{k,m}(s,Z_{s})1(V_{k,m}(s;Z_{s})>\epsilon (1+s) V_{k,m-1}(s;Z_{s}))\right)ds
 \nonumber \\\nonumber \\
=  V_{k,m}(0,Z_0) 
- \mathbb  E\int\limits_0^{\tau_K} c'V_{k,m-1}(s;Z_{s})ds  
 \nonumber \\ \nonumber \\
+ \mathbb  E\int\limits_0^{\tau_K} \frac{k}{1+s}V_{k,m}(s,Z_{s})1(1(1+x_s+y_s > \epsilon (1+s)/k))\, ds. \label{keq}
\end{align}
We estimate the last term here as follows: 
\begin{align*}
\mathbb E\int\limits_0^{\tau_K} \frac{k}{1+s}V_{k,m}(s,Z_{s}) 1(1+x_s+y_s\le \epsilon (1+s)/k)\,ds
 \\\\
\le \mathbb   E\int\limits_0^{\tau_K} \frac{k^{1+b}}{1+s}V_{k,m}(s,Z_{s})\frac{(1+x_s+y_s)^b}{\epsilon^b (1+s)^b}ds
 \\\\
=  k^{1+b}\epsilon^{-b} \mathbb E\int\limits_0^{\tau_K} (1+s)^{k-1-b}V_{m+b}(Z_{s})ds.
\end{align*}
Due to the assumptions, the values $m,b,k$ are such that 
\[
k-1-b < -1, \quad \& \quad m+b \le m_0.
\]
Then, 
\begin{align*}
\mathbb E_{Z_0}\int\limits_0^{\tau_K} (1+s)^{k-1-b}V_{m+b}(Z_{s})ds 
= \mathbb E_{Z_0}\int\limits_0^{\tau_K} (1+s)^{k-1-b}V_{m+b}(Z_{s\wedge \tau_K})ds 
 \\\\
\le  \int\limits_0^{\infty} (1+s)^{k-1-b}\mathbb E_{Z_0}V_{m+b}(Z_{s\wedge \tau_K})ds 
\le  V_{m+b}(Z_0) \int\limits_0^{\infty} (1+s)^{k-1-b}ds = \frac1{b-k} V_{m+b}(Z_0).
\end{align*}
From here and from (\ref{keq}), since $V_{k,m}(0,Z) = V_{m}(Z)$ and since $V_{k,m-1}(s,Z_s)\ge (1+s)^k$ on $(s<\tau_K)$, we conclude (recall that $c' = c-\epsilon = (1-\delta)\gamma - 2m - \epsilon$)
\begin{align}\label{lasteq}
\frac{(1-\delta)\gamma - 2m - \epsilon}{ (k+1)}\mathbb E_{Z_0}\tau_K^{k+1} \le 
\mathbb EV_{k,m}(\tau_K;Z_{t\wedge \tau_K}) 
+ \mathbb E\int\limits_0^{\tau_K} c'V_{k,m-1}(s;Z_{s})ds 
 \nonumber \\\\\nonumber 
\le V_{m}(Z_0) + \frac{k^{1+b}\epsilon^{-b}}{b-k} V_{m+b}(Z_0) 
= V_{m}(Z_0) + \frac{k^{1+m_0-m}\epsilon^{-(m_0-m)}}{m_0-m-k} V_{m_0}(Z_0).
\end{align}
Since $V_{m}(Z_0)\le V_{m_0}(Z_0)$ for $m\le m_0$, the second statement of the Theorem \ref{Thm1} is proved with 
\begin{equation}\label{ckk}
C(k,K) = \frac{(k+1)}{(1-\delta)\gamma - 2m - \epsilon}\left(1+\frac{\epsilon^{-(m_0-m)}k^{1+m_0-m}}{m_0-m-k}\right),
\end{equation}
with any $K, m,m_0,\delta, \epsilon$ satisfying (see (\ref{Kdelta}))
$$
K\ge K(\delta), \;\;
1+k < m < m_0 - k, \;\; (1-\delta)\gamma > 2m_0, \;\;\&\;\; 0< \epsilon < (1-\delta)\gamma - 2m.
$$

~

\noindent
{\bf 3.}
Let $K$ be the large value for which the inequality (\ref{eqtauk}) holds from the part 2 of the Theorem. Of course, it suffices to consider the case $K_1 < K$.  Note that by virtue of the assumptions on the intensities
\begin{equation}\label{q}
q:= \inf_{Z\in {\mathbb K}(K+1)}\mathbb P_Z(\exists \, s\in [0,1]: \, Z_{\tau_K+s}\in {\mathbb K}(K_1)) > 0.
\end{equation}
This is because 
$$
\inf_{Z\in {\mathbb K}(K+2)}(\lambda(Z) \wedge \mu(Z))\ge \frac{\gamma}{1+K+2} >0;
$$
so, the (strong Markov) process which starts at any state in ${\mathbb K}(K+1)$ at any  stopping time (at $T^n$ with any $n$, see below) has a positive probability to hit the set ${\mathbb K}(1)$ over the period of time of length $1/2$, say, and from this set there is again a positive probability to hit the set ${\mathbb K}(K_1)$ over the period of time of length $1/2$ (here we assumed $K_1<1$; if not, then the second transition is not necessary). 
Denote $\Gamma:= {\cal S}\setminus {\mathbb K}(K+1)$ and consider two  sequences of stopping times:
\begin{align*}
\tau^0=T^0:=0, \; \tau^1:=\tau_K, \; T^1:= \inf(t>\tau^1:\; Z_t\in \Gamma) \wedge (\tau^1+1), \; \ldots, 
 \\\\
\tau^{n+1}:=\inf(t>T^{n}:\; Z_t \in {\mathbb K}(K)), \;  T^{n+1}:= \inf(t>\tau^{n+1}:\; Z_t\in \Gamma) \wedge (\tau^{n+1}+1), \; \ldots
\end{align*}
and let
\[
\Delta^n:= \tau^n - T^{n-1}, \;\; n\ge 1.
\]
We have, 
$$
Z_{T^n} \in {\mathbb K}(K+1), \quad \forall \; n, 
$$
and, hence, under $k > 0$,  $\gamma>2m_0 > 2(1+k)$,
\begin{equation}\label{Delta}
\mathbb E_{Z_{T^{n-1}}}(\Delta^n)^{k+1} \le(K+1) C(k,K+1).
\end{equation}
Take any $k<k_1$ so that $2m_0 > 2(1+k_1)$, and denote $\displaystyle p_1=\frac{k_1+1}{k+1}$, and $\displaystyle p_2^{} = \frac{k_1+1}{k_1-k}$.

Due to (\ref{q}), 
\[
\mathbb P(\tau_{K_1}\le T^n | \tau_K> \tau^{n})\ge q>0,
\]
and by induction
\begin{equation}\label{induc}
\mathbb P_Z (\tau^{\ell}<  \tau) \le 
(1-q)^{\ell}.
\end{equation}
So, using strong Markov property, we estimate denoting $\eta_i := \sum_{j=1}^{i}\Delta^j$ by virtue of H\"older's inequality, (\ref{Delta}) and (\ref{induc}),
\begin{align}\label{est1k}
\mathbb E_Z\tau(K_1)^{k+1} 
= \mathbb E_Z\sum_{\ell\ge 1}\tau_{K_1}^{k+1} 1(\tau^{\ell -1} < \tau \le \tau^\ell)
 \nonumber \\ \nonumber \\
\le \sum_{\ell\ge 1}\mathbb  E_Z(\tau_K+(\ell +\eta_\ell))^{k+1}1(\tau^{\ell -1}<  \tau)
  \nonumber \\ \nonumber \\
= \sum_{\ell\ge 1} (\mathbb E_Z (\tau_K+(\ell +\eta_\ell))^{k_1+1})^{1/p_1} (\mathbb E_Z 1(\tau^{\ell -1}<  \tau))^{1/p_2}
  \nonumber \\ \nonumber \\
\le \sum_{\ell\ge 1} ((\ell+1)^{k_1}(\mathbb E_Z \tau_K^{k_1+1}+\ell^{k_1+1} +\eta_\ell)^{k_1+1})^{1/p_1} (\mathbb E_Z 1(\tau^{\ell -1}<  \tau))^{1/p_2}
  \nonumber \\ \nonumber \\
\le  \sum_{\ell\ge 1} [(\ell+1)^{k_1}(C(k_1,K)V_{m_0}(Z) +\ell^{k_1+1} + (K+1)C(k_1,K))]^{1/p_1}(1-q)^{\frac{\ell -1}{p_2}} 
   \nonumber \\ \\ \nonumber 
\le \tilde C(k_1,K) (V_{m_0}(Z)\vee (K+1)), 
\end{align}
as required, e.g., with
$$
\tilde C(k_1,K) \le 
\sum_{\ell\ge 1} [(\ell+1)^{k_1}(C(k_1,K)+\ell^{k_1+1} + C(k_1,K))]^{1/p_1}(1-q)^{\frac{\ell -1}{p_2}}.
$$
All the statements of the Theorem are thus  proved.

\medskip

\begin{remark}
It is tempting to claim that the random variables  $\eta_\ell$ and $1(\tau^{\ell -1}<  \tau)$ are independent. If this were correct, then the last estimate would have been much better (with a much smaller value of the constant) and the last calculus a bit easier; also the analogue of the inequality (\ref{eqtauk1}) for $k=0$ would have been valid under the assumption $\gamma > 2m_0 \ge 2$. Yet, in our model we do not see how to justify this fairly  plausible  claim. 
\end{remark}

\subsection*{Note on convergence rate}
Now assume that we achieved the situation that two independent strong Markov processes $Z$ and $Z'$  both  attain the set $\mathbb K$ at some stopping time $T=T_1$. Then, due to the assumption (\ref{usl}) by changing appropriately the probability space (which must be explained in the full presentation) we manage to arrange gluing the two equivalent processes with a positive probability bounded away from zero on the interval $[T,T+1]$; if at least one of the processes leaves the compact ${\mathbb K}(K_1+1)$, we stop the couple at the moment of exit. If at this step coupling was not successful, we wait till they both attain $\mathbb K$ for the next time $T_2$ , etc. In this way and using the analogue of the moment inequality for the pair of two independent copies of our Markov process one of which is stationary, it is possible to establish a polynomial bound for  the convergence rate towards the stationary distribution (which again shows that it is unique), 
\[
\|{\cal P}_t^Z - \pi\|_{TV} \le C(Z, \delta)(1+t)^{-k+\delta} 
\]
with any $\delta>0$, 
where the norm in the left hand side is in total variation, and where ${\cal P}_t^Z$ is the distribution of the process given the initial state $Z$, and $\pi$ is a stationary distribution of the process. 
We postpone it till further studies.

\subsection*{Acknowledgement} 
The author is grateful to the anonymous referee  for very helpful comments.

\end{document}